\documentclass[a4paper,fleqn]{cas-sc}

\usepackage[authoryear]{natbib}
\usepackage{tikz}
\usetikzlibrary{positioning,arrows.meta,shapes.geometric}
\usepackage{algorithm,algpseudocode}
\usepackage{float}
\usepackage{graphicx}
\usepackage{orcidlink}

\def\tsc#1{\csdef{#1}{\textsc{\lowercase{#1}}\xspace}}
\tsc{WGM}
\tsc{QE}
\tsc{EP}
\tsc{PMS}
\tsc{BEC}
\tsc{DE}
\begin{document}

\let\WriteBookmarks\relax
%\def\floatpagepagefraction{1}
%\def\textpagefraction{.001}

\shorttitle{Arctic Oil Spill Response Model}

\shortauthors{Rahman et~al.}

\title [mode = title]{A Two-Stage Stochastic Optimization Framework for Environmentally Sensitive Oil Spill Response Resource Allocation in the Arctic}                      
%\tnotemark[1,2]

%\tnotetext[1]{This document is the results of the research
   %project funded by the National Science Foundation.}

%\tnotetext[2]{The second title footnote which is a longer text matter
   %to fill through the whole text width and overflow into
   %another line in the footnotes area of the first page.}

\makeatletter
\@ifundefined{frontmatter@orcid}{}{\renewcommand*\frontmatter@orcid{}}
\@ifundefined{printorcid}{}{\renewcommand\printorcid[1]{}}
\@ifundefined{printOrcid}{}{\renewcommand\printOrcid[1]{}}
\@ifundefined{orcid}{}{\renewcommand\orcid[1]{}}
\makeatother


\author[]{Md Ashiqur Rahman\,\orcidlink{0000-0001-5135-4149}}
\ead{ashiqur.txst@gmail.com}
\author{Mustofa Tanbir Kuhel\,\orcidlink{0000-0001-5734-7679}}
\author{Clara Novoa\,\orcidlink{0000-0001-8502-5335}}
\ead{cn17@txstate.edu}

\affiliation{organization={Ingram School of Engineering, Texas State University},
  addressline={601 University Drive}, city={San Marcos}, state={Texas},
  postcode={78666}, country={USA}}

\cortext[1]{Corresponding author} % single line, no \\

%\cormark[1]


%\ead{cvr_1@tug.org.in}

%\ead[url]{www.cvr.cc, cvr@sayahna.org}

%\credit{Conceptualization of this study, Methodology, Software}

%\ead{cvr3@sayahna.org}
%\ead[URL]{www.sayahna.org}

%\credit{Data curation, Writing - Original draft preparation}

% % Fourth author
% \author%
% [1,3]
% {Rishi T.}
% \cormark[2]
% \fnmark[1,3]
% \ead{rishi@stmdocs.in}
% \ead[URL]{www.stmdocs.in}

% \affiliation[3]{organization={STM Document Engineering Pvt Ltd.},
%     addressline={Mepukada}, 
%     city={Malayinkil},
%     % citysep={}, % Uncomment if no comma needed between city and postcode
%     postcode={695571}, 
%     state={Trivandrum},
%     country={India}}


%\fntext[fn2]{Another author footnote, this is a very long footnote and
  %it should be a really long footnote. But this footnote is not yet
  %sufficiently long enough to make two lines of footnote text.}

%\nonumnote{This note has no numbers. In this work we demonstrate $a_b$
  %the formation Y\_1 of a new type of polariton on the interface
  %between a cuprous oxide slab and a polystyrene micro-sphere placed
  %on the slab.
  %}

\makeatletter
\let\oldprintFirstPageNotes\printFirstPageNotes
\renewcommand{\printFirstPageNotes}{%
  \vspace*{-6pt}% reduce space before the block
  \oldprintFirstPageNotes
  \vspace*{-6pt}% reduce space after the block
}
\makeatother


\begin{abstract}
The risk of oil spills in the Alaskan Arctic has become an urgent environmental and logistical concern as maritime traffic increases due to climate-induced sea-ice retreat. Traditional deterministic response-planning models fail to capture the region’s complex uncertainties, including variable spill magnitudes, fluctuating environmental sensitivities, and infrastructure limitations. This study develops a two-stage stochastic mixed-integer linear programming framework to optimize the location of oil-spill response stations and the allocation of heterogeneous resources across multiple probabilistic spill scenarios. The model integrates weighted objectives combining spill volume, environmental sensitivity index (ESI), response time, and operational costs related to station setup, deployment, and inter-station transfers. Distinct importance weight and ecological weight balance ecological protection and operational efficiency. Data were compiled from Alaska Department of Environmental Conservation (ADEC) spill records and National Oceanic and Atmospheric Administration (NOAA) ESI layers and engineered into model-ready scenarios through data harmonization and sampling. The model was solved using the Gurobi optimizer, and sensitivity analysis was conducted across 324 configurations of importance weight ($k_1$,$ k_2$), and internal ecological weight ($\omega_1$, $\omega_2$, $\omega_3$). Results show a 35.45 \% improvement in response effectiveness over deterministic methods which is verified by the Value of the Stochastic Solution (VSS), demonstrating the model’s capacity to prioritize high risk and sensitive spill events. Pareto frontier revealed actionable trade-offs between cost efficiency and ecological coverage, with environmental sensitivity ($\omega_2\ = 0.8$) emerging as the dominant factor in optimal result. The proposed framework equips Arctic emergency planners with a data-driven decision-support tool that balances ecological protection, and logistical feasibility under uncertainty.
\end{abstract}

% \includegraphics{figs/grabs.pdf}

%\item Research highlights item 1
%\item Research highlights item 2
%\item Research highlights item 3

% Keywords
\begin{keywords}
Oil Spill Response \sep Environmental Sensitivity \sep Response Time \sep Two-stage Stochastic Programming \sep Mixed-Integer Linear Programming \sep Decision Support Tool
\end{keywords}

\maketitle

\section{Introduction}
\vspace{0.5em}
Maritime activity in the Alaskan Arctic has increased substantially over the past decade as declining sea ice has expanded operational windows and enabled greater commercial access. Vessel transits through the Bering Strait more than doubled between 2009 and 2021, rising from 262 to over 545 annual crossings (\cite{GAO2022}; \cite{WWF2022CrossingTheLine}). This growth elevates the probability that oil spills will occur near ecologically fragile shorelines and regions that support Indigenous subsistence practices. Such risks are amplified by the Arctic’s unique environmental and infrastructural constraints, which make effective response operations inherently difficult.\\

Oil spill response in northern Alaska is characterized by long mobilization distances, limited staging infrastructure, severe weather, and seasonally variable sea ice. These constraints slow the deployment of critical equipment and extend the duration of contamination, with disproportionate consequences for coastal communities dependent on marine resources (\cite{ArcticCouncil_UnknownTitle}; \cite{Afenyo2020RiskArcticOilSpills}). The grounding of the Kulluk drilling rig in 2012 illustrates these challenges: although the incident occurred within proximity of a Coast Guard station, adverse conditions and limited logistical capacity delayed response efforts for up to four days (\cite{TSAC2014Task14-01}; \cite{NTSB2015MAB-15-10}). Broader assessments similarly highlight persistent gaps in governance, infrastructure, and operational readiness that constrain Alaska’s ability to respond to large-scale emergencies (\cite{Daisy2022CanadianMarineOilSpill}; \cite{ChircopLEsperance2016FunctionalInteractionsMaritimeRegulation}).\\

Mitigating these risks requires planning frameworks that explicitly incorporate the uncertainty inherent in spill occurrence, environmental conditions, and accessibility. At the same time, response strategies must reflect the environmental sensitivity of Arctic coastlines, as documented through NOAA Environmental Sensitivity Index (ESI) maps, which classify shorelines, habitats, and human-use resources based on their vulnerability to oil contamination (\cite{DAffonseca2023EnvironmentalSensitivityIndex}; \cite{NOAAESI2016}). Operational guidance from Alaska Clean Seas further outlines realistic deployment capabilities, equipment performance, and setup requirements (\cite{AlaskaCleanSeas2021TacticsDescriptions}). However, these data sources are rarely integrated within unified decision-support models for the region.\\

Given these considerations, a key limitation in existing planning methodologies becomes evident. There is no model that concurrently incorporates spill uncertainty, region-specific logistical constraints, and ESI-based ecological priorities for the Alaskan Arctic. Deterministic siting models provide coarse baselines but underestimate variability. Existing stochastic models handle uncertainty but treat environmental sensitivity externally or simplify operational constraints. Consequently, decision-makers lack a comprehensive tool that captures the full interplay of uncertainty, logistics, and ecological consequences.\\

This study addresses this gap by developing a two-stage stochastic mixed-integer programming model tailored to oil spill response planning in northern Alaska. The model determines station locations and resource prepositioning decisions in the first stage, and in the second stage allocates and transfers equipment across probabilistic spill scenarios that vary in spill size, environmental sensitivity, and response time requirements. By integrating realistic logistics, sensitivity-driven priorities, and stochastic spill behavior, the framework provides a more robust basis for evaluating response readiness under Arctic conditions.\\

The remainder of this paper is organized as follows. \textbf{Section 2} reviews the literature on deterministic and stochastic facility-location and oil spill response models, with emphasis on Arctic applications and the treatment of environmental sensitivity. \textbf{Section 3} presents the two-stage stochastic mixed-integer linear programming model, introduces the notation, and describes the construction of spill scenarios and associated data sources. \textbf{Section 4} reports the base-case optimization results, including station selection, multi-resource allocation patterns, and value-of-information metrics such as the Value of the Stochastic Solution (VSS) and the Expected Value of Perfect Information (EVPI). \textbf{Section 5} conducts sensitivity analysis on global coverage–cost weights and internal ecological weights, and develops the Pareto frontier between ecological coverage and cost. 

\section{Literature Review}
\vspace{0.5em}

Research on emergency facility siting and response logistics has traditionally relied on deterministic optimization models. \cite{ChurchVelle1974MaximalCovering} introducted the Maximal Covering Location Problem (MCLP) that forms a foundational framework for locating facilities to maximize service coverage. In the context of oil spill response, several studies extend the maximal-covering and partial-covering paradigm to operational decision support. \cite{BelardoHarraldWallaceWard1984PartialCovering} introduced one of the earliest mixed-integer formulations for allocating maritime response assets under partial-coverage criteria. Building on this systems-oriented perspective, \cite{KeramitsoglouCartalisKassomenos2003DSSOilSpill} developed a decision-support system that integrates equipment positioning with real-time environmental inputs. More recent work, such as \cite{Medini2018SpatialDSSOilSpill}, incorporates shoreline sensitivity mapping into a deterministic siting structure, thereby linking environmental vulnerability with strategic resource deployment. Although these contributions remain influential, they assume fixed spill locations, constant travel times, and full availability of response assets and these assumptions are unrealistic under Arctic operating conditions.\\

Scenario-based and stochastic approaches have been introduced to overcome deterministic limitations. \cite{DordevicSabaliaMohovicBrcic2022SkimmerSelection} modeled skimmer deployment under uncertain spill trajectories. \cite{AmirHeidariRaie2019StochasticOilSpillRiskPersianGulf} developed a multi-scenario allocation framework for the Persian Gulf. \cite{GarrettSharkeyGrabowskiWallace2017DynamicResourceAllocationOilSpill} proposed a stochastic resource allocation model for Arctic spill response that accounts for variable travel times and environmental conditions. Reviews by \cite{ShuCuiSongGanXuWuZheng2024InfluenceSeaIce} emphasize the necessity of embedding uncertainty, ice dynamics, and harsh weather into Arctic response models. Despite these developments, ecological sensitivity considerations are often handled outside the optimization framework or represented in simplified form.\\

Environmental sensitivity mapping, particularly through NOAA’s ESI program, provides structured data on vulnerable shoreline types, species habitats, and socio-economic resources (\cite{DAffonseca2023EnvironmentalSensitivityIndex}; \cite{NOAAESI2016}). Complementary classifications exist in the Canadian Arctic (\cite{DFO2025ArcticRegion}). While these datasets are widely used in environmental planning, they are not routinely integrated into optimization-based response models. \cite{DasGoerlandtPelot2025FacilityLocationArcticOilSpill} used sensitivity weights in a Bayesian decision model, but their operational framework remained deterministic. \cite{DasGoerlandtPelot2024MixedIntegerProgrammingOilSpillArctic} introduced a mixed-integer model that jointly considers spill size, ESI-based sensitivity, and response time in the Canadian Arctic which demonstrates meaningful performance improvements over classical MCLP baselines. However, this model does not incorporate Alaska’s unique logistical barriers such as limited port infrastructure, dependence on airlift, and discontinuous road networks.\\

Additional studies underscore the specific challenges of the Alaskan Arctic. Government assessments document major infrastructure deficits, operational fragility, and geographically constrained response capacity. The Kulluk incident in 2012 serves as a clear case study of how severe conditions can delay response even near established assets. Reports from the World Wildlife Fund highlight rapidly increasing shipping pressure and its ecological implications for the Bering Strait corridor. Operational manuals from Alaska Clean Seas detail equipment limitations, setup times, and deployment strategies that are largely absent from optimization studies.\\

Taking together, this body of literature reveals a systematic gap. Existing models address either uncertainty, environmental sensitivity, or Arctic-specific logistics, but no study integrates all three within a single optimization framework for the Alaskan Arctic. Deterministic approaches fail to capture variability while stochastic models lack explicit ecological prioritization, and sensitivity-focused studies do not incorporate realistic operational constraints. The present study fills this gap by combining: (i) scenario-based uncertainty in spill size, timing, and location, (ii) Environmental Sensitivity Index–driven ecological prioritization, and (iii) logistically realistic deployment and transfer capabilities specific to northern Alaska.

\section{Method and Data}
\vspace{0.5em}

The problem most closely resembles  Stochastic Maximal Covering Location Problem (SMCLP), which is an extension of the classical Maximal Covering Location Problem (MCLP) (\cite{Daskin1983}). This study employs a two-stage stochastic mixed-integer linear programming (MILP) model to address the dual challenges of uncertainty and ecological sensitivity in Arctic oil spill response. The first stage determines the optimal selection of 10 response stations from candidate locations (e.g., station1–station10) which incorporates fixed setup costs (\$1,300–\$3,700) and geospatial accessibility. The second stage allocates resources dynamically across three spill scenarios (low, medium, high severity), each weighted by probability (20\%, 50\%, 30\%) and environmental sensitivity indices (ESIs).

\subsection{Problem Statement}
\vspace{0.5em}
Traditional deterministic models, which assume fixed spill magnitudes and static resource demands, cannot adequately represent Arctic oil spill conditions where spill size, ice state, and environmental sensitivity are highly variable. Existing deterministic frameworks either ignore uncertainty or treat it through simple safety factors, and sensitivity-focused approaches prioritize ecologically vulnerable areas without modeling probabilistic spill occurrence. Stochastic models proposed in the broader oil spill literature address scenario uncertainty but generally neglect Arctic-specific constraints such as sparse port infrastructure, long mobilization distances, and limited road connectivity. As a result, decision-makers lack a tool that simultaneously captures spill uncertainty, spatially explicit ecological priorities, and realistic Arctic logistics when planning response station locations and resource allocations.\\

This study therefore formulates a decision problem in which a planner must determine which response stations to open and how to allocate heterogeneous resources so that the expected weighted coverage of oil spills is maximized under uncertainty. The model must balance three competing criteria: spill volume, environmental sensitivity, and response time, while respecting constraints on station budgets, resource capacities, inter-station transfers, and maximum allowable arrival times. The subsequent subsections present a two-stage stochastic mixed-integer linear programming formulation of this problem, including the definition of decision variables, parameters, and governing constraints.

\begin{figure}
    \centering
    \includegraphics[width=0.5\linewidth]{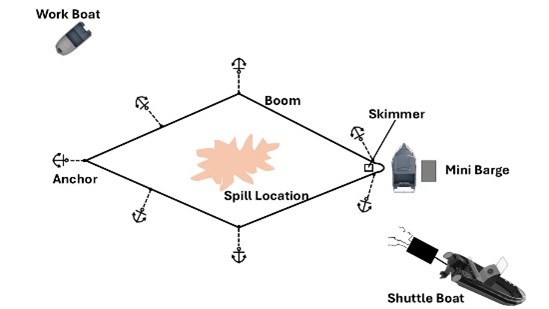}
    \caption{Illustration of oil spill recovery}
    \label{fig:pf1}
\end{figure}

\subsection{Location-Allocation Model}
\vspace{0.5em}
The proposed model is formulated as a two-stage stochastic mixed integer linear programming problem. In the first stage, the model determines the optimal locations of oil spill response stations under a predefined station budget. In the second stage, it allocates heterogeneous response resources like skimmers, booms, and dispersants to spill sites across multiple scenarios characterized by varying spill sizes, environmental sensitivity indices, and response time requirements. The objective function maximizes expected spill coverage which is weighted by spill volume, ecological sensitivity, and timeliness while minimizing total costs from station setup, resource deployment, and inter-station transfers. A set of operational constraints ensures feasible station activation, resource capacity limits, transfer logistics, and adherence to maximum allowable response times. The Arctic oil spill response model is formulated using a structured set of elements that integrates the spatial, operational, and stochastic complexities existing in emergency resource planning.  Figure \ref{fig:pf1} provides an illustrative representation of the oil spill recovery approach. The model notation and terminologies are on the table below. A pairing set $P \subseteq O \times I$ specifies which stations are eligible to respond to each spill based on logistical feasibility and geographic proximity. 

\subsubsection{Mathematical Model}

\begin{align}
\max\; Z
=\;& k_1 
      \sum_{k \in K} \rho^{k} 
      \sum_{i \in I} \sum_{o \in O}
      \left( 
        \omega_1 v_{o}^{k} 
        + \omega_2 \eta_{o}^{k} 
        - \omega_3 T_{o}^{k} 
      \right)
      Y_{io}^{k} \notag \\
&\; - k_2 
  \left[
    \sum_{k \in K} \rho^{k}
    \left(
         \sum_{r \in R} \sum_{i \in I} \sum_{o \in O}
                dc_{io}\, Z_{io}^{kr}
       +
         \sum_{i \in I} \sum_{j \in I} \sum_{r \in R}
                tc_{ij}\, A_{ij}^{kr}
    \right)
    + \sum_{i \in I} C_i X_i
  \right]
\tag{1}
\end{align}

\begin{align}
X_i &\ge Y_{io}^{k},
& & \forall (o,i)\in P,\ \forall k\in K 
\tag{2}
\\[2pt]
\sum_{i\in I} X_i &\le n,
\tag{3}
\\[2pt]
\sum_{i\in I} Y_{io}^{k} &\le 1,
& & \forall o\in O,\ \forall k\in K
\tag{4}
\\[4pt]
Z_{io}^{kr} &\le d_{o}^{kr} Y_{io}^{k},
& & \forall i\in I,\ \forall o\in O,\ \forall k\in K,\ \forall r\in R
\tag{5}
\\[4pt]
Z_{io}^{kr} &\le R_i^{r} Y_{io}^{k},
& & \forall i\in I,\ \forall o\in O,\ \forall k\in K,\ \forall r\in R
\tag{6}
\\[4pt]
\sum_{o\in O} Z_{io}^{kr}
&\le R_i^{r}
   + \left(
        \sum_{j\in I} A_{ij}^{kr}
        - \sum_{j\in I} A_{ji}^{kr}
     \right) X_i,
& & \forall i\in I,\ \forall k\in K,\ \forall r\in R
\tag{7}
\\[2pt]
A_{ij}^{kr} &\le R_i^{r} X_i,
& & \forall i\in I,\ \forall j\in I,\ \forall k\in K,\ \forall r\in R
\tag{8}
\\[4pt]
\sum_{i\in I} Z_{io}^{kr} 
&\ge d_{o}^{kr} \sum_{i\in I} Y_{io}^{k},
& & \forall o\in O,\ \forall k\in K,\ \forall r\in R
\tag{9}
\\[4pt]
pt_i^{k} 
&\le d + \sum_{r\in R} d_{o}^{kr} st_r,
& & \forall o\in O,\ \forall k\in K,\ \forall i\in I
\tag{10}
\\[4pt]
T_{o}^{k} 
&\ge (\theta_{io} + pt_i^{k}) Y_{io}^{k},
& & \forall o\in O,\ \forall k\in K,\ \forall i\in I
\tag{11}
\\[4pt]
T_{o}^{k} 
&\le \tau_{\max} \sum_{i\in I} Y_{io}^{k},
& & \forall o\in O,\ \forall k\in K
\tag{12}
\\[4pt]
X_i &\in \{0,1\},\quad
Y_{io}^{k} \in \{0,1\},\quad
Z_{io}^{kr} \ge 0,\quad
A_{ij}^{kr} \ge 0,\quad
T_{o}^{k} \ge 0,
& &
\tag{13}
\end{align}

\begin{table}[width=.9\linewidth,cols=2,pos=htbp]
	\caption{Sets, parameters, and decision variables of the stochastic oil spill response model.}
	\label{tbl:sets_params}
	\begin{tabular*}{\tblwidth}{@{} LL@{}}
		\toprule
		\textbf{Notation} & \textbf{Description} \\
		\midrule
		\multicolumn{2}{@{}l}{\textit{Sets}} \\[2pt]
		$O$ & Hypothetical oil spills, indexed by $o$ \\[2pt]
		$I$ & Potential response stations, indexed by $i$ \\[2pt]
		$K$ & Scenarios, indexed by $k$ \\[2pt]
		$P$ & Station--spill index pairs $(o,i)\in O\times I$ \\[2pt]
		$R$ & Resource types, indexed by $r$ \\
		\midrule
		\multicolumn{2}{@{}l}{\textit{Parameters}} \\[2pt]
		$R_i^{r}$ & Total available quantity of resource $r$ at station $i$ \\[2pt]
		$\rho^{k}$ & Probability of scenario $k$ \\[2pt]
		$\eta_{o}^{k}$ & Environmental sensitivity index of spill $o$ in scenario $k$ \\[2pt]
		$\tau_{\max}$ & Maximum allowable response time \\[2pt]
		$\omega_1,\omega_2,\omega_3$ & Weights for spill size, environmental sensitivity, and response time \\[2pt]
		$k_1,k_2$ & Relative importance factors for coverage and cost \\[2pt]
		$n$ & Number of response stations to be selected \\[2pt]
		$d$ & Base delay time (crew readiness, administrative delay, etc.) \\[2pt]
		$st_{r}$ & Setup time for resource type $r$ \\[2pt]
		$v_{o}^{k}$ & Size (volume) of spill $o$ in scenario $k$ \\[2pt]
		$C_i$ & Opening cost of station $i$ \\[2pt]
		$d_{o}^{kr}$ & Required quantity of resource $r$ to cover spill $o$ in scenario $k$ \\[2pt]
		$dc_{io}$ & Deployment cost per unit of resource from station $i$ to spill $o$ \\[2pt]
		$tc_{ij}$ & Transfer cost per unit of resource from station $i$ to station $j$ \\[2pt]
		$pt_i$ & Preparation time for station $i$ (hours) \\[2pt]
		$\theta_{io}$ & Travel time from station $i$ to spill site $o$ \\[2pt]
		\midrule
		\multicolumn{2}{@{}l}{\textit{Decision variables}} \\[2pt]
		$X_i$ & $1$ if response station $i$ is selected\\[2pt]
		&  $0$ otherwise\\[2pt]
		$Y_{io}^{k}$ & $1$ if station $i$ is selected to cover spill $o$ in scenario $k$\\[2pt]
		&  $0$ otherwise\\[2pt]
		$Z_{io}^{kr} \in R_{0}^+$ & Amount of resource type $r$ deployed from station $i$ to spill $o$ in scenario $k$ \\[2pt]
		$A_{ij}^{kr} \in R_{0}^+$ & Amount of resource type $r$ transferred from station $i$ to station $j$ in scenario $k$ \\[2pt]
		$T_{o}^{k} \in R_{0}^+$ & Arrival (response) time at spill site $o$ in scenario $k$ \\
		\bottomrule
	\end{tabular*}
\end{table}

Expression (1) denotes the objective function of the model. Constraint (2)  ensures spills are only assigned to active stations.  Constraint (3)  limits the number of open facilities. Constraint (4) restricts each spill to be covered by at most one station. Constraint (5) ensures deployment demand limit. Constraint (6) enforces that dispatched resources do not exceed either site specific demand or station capacity respectively. Constraint (7) regulates resource movement. Constraint (8) ensures station inventories remain feasible. Constraint (9) ensures that the least required amount of resources is deployed to spill sites. Constraint (10) accounts for operational delays. Constraint (11) keeps track of arrival time of response team at spill site. Constraint (12) response time window constraints enforce maximum allowable arrival times. Finally, constraint (13) defines the binary or continuous nature of all decision variables. 

\subsection{Data}
\vspace{0.5em}
This study utilizes real-world spill incident records obtained from the Alaska Department of Environmental Conservation (ADEC) Spill Search database (\cite{alaskaSpillSearch}). 

\subsubsection{Stations and resources}
\vspace{0.5em}
The strategic placement of response stations in the Arctic is informed by geographic accessibility, proximity to high-risk shipping corridors, and ecological vulnerability. Studies, such as those by the Alaska Department of Environmental Conservation (ADEC). There are various types of spills. This paper focuses on the spill incidents for oil and that recovered by the response station by field visit. \cite{Etkin2004} proposed a model to find the response cost using Basic Oil Spill Cost Estimation Model (BOSCEM). Although the total number of resources such as Skimmer, Boom, Dispersant in every station is confidential, and not easily accessible, according to Alaska Clean Seas (ACS) catalogue, currently, the value of dedicated spill response equipment warehoused under ACS management exceeds \$75 million. ACS owns around half of this equipment, with the remaining half owned by member companies such as Alpine, Kuparuk, and Marine Base. \\

Using the per-unit resource price and total station cost, the number of resources in each station is also assumed (\cite{Morris2011TRB}). In this study, interstation transfer costs represent the operational expenses incurred when transporting oil spill response equipment such as skimmers, booms, and dispersants between selected response stations. Transfer cost, $tc_{ij}$ were derived using a distance-based approach where Euclidean distances (in kilometers) between stations were first computed to construct a symmetric distance matrix. A uniform multiplier of $\alpha =1.723$ was applied to convert these distances into monetary transfer cost to each non-diagonal entry which yields the transfer cost formula $tc_{ij}=1.723 d_{ij}$ for all $i \neq j$. This scalar was empirically derived by scaling the original distance matrix to match realistic transfer costs and is consistent with reported land transportation rates in industry and federal guidance. For example, the Western Canada Marine Response Corporation (WCMRC) reports truck-based equipment transport at \$0.65–\$1.30/km depending on vehicle class and load (\cite{WCMRC2024}), while U.S. General Services Administration mileage rates and Transport Research Board (TRB) reports cite similar or higher rates for hazardous and oversized cargo up to \$2.57/km for dispersant tankers (\cite{GSA_POV_Rates}). \\

According to standard operating ratios mentioned in oil spill literature and guideline documents, the usual dispersant-to-oil application ratio varies from 1:10 to 1:50, depending on oil type, weathering, and environmental circumstances (\cite{Fingas2017_OilSpillScienceTechnology}; \cite{NOAA_Marshes_OilSpills}). For modeling consistency and conservative estimation, this study uses a fixed dispersant application ratio of 1:50. This study assumed 10 feet of boom per gallon of oil spilled. A single skimmer unit can efficiently recover up to several hundred gallons per operational shift, depending on oil kind and sea state (\cite{EPA1999_OilSpillResponse}). Table \ref{tbl:stations} presents the station name and amount of resources estimated.

\begin{table}[width=.9\linewidth,cols=6,pos=h]
\caption{Response station locations and available resource inventories.}
\label{tbl:stations}
\begin{tabular*}{\tblwidth}{@{} L L C C C C @{}}
\toprule
\textbf{Station Name} & \textbf{Location} & \textbf{Cost (M\$)} & 
\textbf{Boom (ft)} & \textbf{Skimmers} & \textbf{Dispersants (gal)} \\
\midrule
Deadhorse 
& 70.1952$^\circ$N, $-$148.4651$^\circ$W 
& 41.25 
& 200{,}000 
& 30 
& 7{,}000 
\\[4pt]

A4W1 Marine Base 
& 70.2508$^\circ$N, $-$148.4696$^\circ$W 
& 15 
& 100{,}000 
& 20 
& 12{,}000 
\\[4pt]

Kuparuk Field Site 
& 70.3300$^\circ$N, $-$149.6100$^\circ$W 
& 11.25 
& 80{,}000 
& 15 
& 5{,}000 
\\[4pt]

Alpine Field Site 
& 70.3420$^\circ$N, $-$150.9470$^\circ$W 
& 7.5 
& 40{,}000 
& 8 
& 2{,}000 
\\
\bottomrule
\end{tabular*}
\end{table}

\subsubsection{Spill Events}
\vspace{0.5em}
Spill location and spill size data were given on the dataset. Appropriate spill location and spill size data were extracted with respective to their spill id. For modeling oil spill containment in Arctic waters, a standardized boat speed of 1.85 km/h (1 knot) is adopted. For each station, the total preparation time is computed as follows:

\begin{equation}
\text{prep\_time}_{i}
    = \frac{\text{base\_delay} 
      + \sum_{r \in R} \text{units}_{or} \cdot \text{setup\_time}_{r}}
      {60}
\tag{14}
\end{equation}

Where $unit_{or}$  denotes the number of units of resource $r$ required by spill $o$, $setup\_time_r$   is the estimated time in minutes required to prepare a single unit of resource $r$, and  $base\_delay_s$ represents a fixed mobilization delay at station $s$ accounting for personnel readiness and administrative tasks. We assume this base delay time as 15 minutes base delays, and the setup is based on the citation (\cite{ArcticCouncil_UnknownYear}).\\

The result is expressed in hours to ensure consistency with the travel time component in the total response time model.

\subsubsection{Scenario Generation}
\vspace{0.5em}
A hybrid scenario generation framework was developed using Monte Carlo simulation to capture the embedded uncertainty in Arctic oil spill response planning. This framework produced five scenarios: one deterministic scenario based on original historical observations, and four stochastic scenarios derived through sampling from statistically fitted and empirically corrected distributions. Distribution fitting was conducted using the distfit Python package with model selection guided by the Kolmogorov–Smirnov (KS) test for goodness of fit. For spill size and boom demand, exponential distributions were identified as the best fit with acceptable statistical performance. The statistics are as follows:

\begin{table}[width=.6\linewidth,cols=4,pos=h]
\caption{Scenario generation}
\label{tbl:scenario_generation}
\begin{tabular*}{\tblwidth}{@{} L C C C C @{}}
\toprule
KS statistic & p-value & spills & scenarios & scenario instances \\
\midrule
0.1126 & 0.3725 & 17 & 5 & 85 \\
\bottomrule
\end{tabular*}
\end{table}

Scenario 1 retained the original deterministic values, while Scenarios 2–5 represented independent stochastic realizations showing the variability of key parameters. Validation of the stochastic scenarios demonstrated good agreement with original data with mean and standard deviation deviates generally within 2–15\% for most parameters. The resulting scenario set supports stochastic optimization and enables subsequent evaluation of the VSS to quantify the benefit of explicitly modeling uncertainty in oil spill response decision-making.

\subsection{Solution Approach}
\vspace{0.5em}
The solution procedure integrates data validation, scenario generation, and a two–stage stochastic mixed-integer linear programming (MILP) model that jointly determines station selection, resource deployment, and inter-station transfers under uncertainty. The corresponding computational steps are formalized in Algorithm~\ref{alg:stochastic_spill_response}.

The process begins with the ingestion of all available operational and environmental data, including potential station locations, resource inventories, spill characteristics, and travel time parameters. After that, a data sufficiency check is performed to ensure that all essential information required for model construction is available. If deficiencies are detected, the procedure iteratively acquires or updates the missing data and revalidates the dataset before proceeding. Once the dataset is deemed sufficient, a set of stochastic scenarios is generated. Each scenario $k \in K$ represents a possible realization of spill size $v_o^k$, environmental sensitivity $\eta_o^k$, required resource demand $d_o^{kr}$, and travel-time variations $\theta_{io}$. The sampling process captures the complex uncertainty inherent in Arctic spill events and assigns an occurrence probability $\rho^k$ to each scenario.

The core of the solution approach is an iterative optimization loop. In each iteration, the two stage stochastic MILP is constructed using the previously generated scenario set. The first-stage variables $X_i$ determine which response stations are opened, while the second-stage recourse variables $(Y_{io}^k, Z_{io}^{kr}, A_{ij}^{kr}, T_o^k)$ capture spill–station assignments, resource deployments, inter-station transfers, and arrival times for every scenario. The objective maximizes a weighted coverage metric while penalizing deployment, transfer, and station establishment costs. Constraints ensure feasible resource allocation, station capacity compliance, inventory balance, and adherence to response-time limits.

The model is solved using a commercial MILP solver (Gurobi) which produces candidate solution for station locations and resource allocations. The obtained solution is then evaluated for inventory feasibility, demand satisfaction, and time-window compliance. If any station violates capacity, or if any spill exceeds allowable response-time limits, a corrective mechanism is invoked. This mechanism modifies the candidate station set, for example, by activating
additional stations, adjusting allowable selections, or expanding inventory at critical locations, and the two stage stochastic MILP is resolved. This loop continues until all stations satisfy coverage and resource sufficiency conditions across all scenarios.

Once feasibility is achieved, the algorithm returns the optimal station configuration and resource allocation plan. The final outputs include spill coverage statistics, scenario-specific deployment plans, and performance metrics such as RP, EWS, EVP, EEV, VSS, and EVPI, which quantify the benefit of using the stochastic model relative to deterministic alternatives. These outputs provide actionable insights for emergency planners by revealing the trade-offs between cost, coverage, and uncertainty, and by identifying conditions under which additional information or strategically chosen station locations offer the greatest marginal value.

\begin{algorithm}[h]
\caption{Stochastic Oil Spill Response Procedure}
\label{alg:stochastic_spill_response}
\begin{algorithmic}[1]
\Require Sets $O,I,R,K$; parameters $\{v_o^k,\eta_o^k,d_o^{kr},\theta_{io},
        R_i^r,dc_{io},tc_{ij},C_i,\rho^k\}$  
\Ensure $(X_i^\ast,Y_{io}^{k\ast},Z_{io}^{kr\ast},A_{ij}^{kr\ast},T_o^{k\ast})$

\State \textbf{Input} all data $\mathcal{D}$
\While{\textsc{DataCheck}$(\mathcal{D}) = \text{false}$}  
    \State $\mathcal{D} \leftarrow$ \textsc{UpdateData}$(\mathcal{D})$
\EndWhile

\For{$k \in K$}
    \For{$o \in O$}
        \State $(v_o^k,\eta_o^k,d_o^{kr},\theta_{io}) \sim \mathcal{S}$ 
    \EndFor
\EndFor

\State $\text{feasible} \leftarrow \text{false}$

\While{\text{feasible} = \text{false}}
    \State \textbf{Build} model
        \Statex \hspace{1.6em} $\max Z(X,Y,Z,A,T)$ given eqn. (13) - (25)
    \State $(X^\ast,Y^\ast,Z^\ast,A^\ast,T^\ast) \leftarrow
           \textsc{SolveMILP}(O,I,R,K)$
    \State $(\Delta R_i^r,\Delta d_o^{kr},\Delta T_o^k)
           \leftarrow \textsc{Evaluate}(X^\ast,Y^\ast,Z^\ast,A^\ast,T^\ast)$
    \If{$\Delta R_i^r \le 0,\ \Delta d_o^{kr} \le 0,\ T_o^k \le \tau_{\max},
         \ \forall i,o,r,k$}
        \State $\text{feasible} \leftarrow \text{true}$
    \Else
        \State $I \leftarrow \textsc{AdjustStations}(I,X^\ast,\Delta R)$
    \EndIf
\EndWhile

\State \textbf{Output} $(X^\ast,Y^\ast,Z^\ast,A^\ast,T^\ast)$ and Evaluation metrics 
      (RP, EWS, EVP, EEV, VSS, EVPI)

\end{algorithmic}
\end{algorithm}

\section{Results}
\vspace{0.5em}
The optimization model selected station 3 and station 4 as the optimal response stations to be established. These stations were chosen from available four stations which are based on their strategic proximity to high-sensitivity spill zones, existing resource inventories, and network centrality. The total fixed cost associated with establishing these two stations was \$18.75 million, with station3 bearing the larger portion (\$11.25 million).\\

Across all five scenarios, the selected stations provided full resource coverage for 30 of the 85 possible spill-scenario combinations which result in an overall coverage rate of 35.3\%. Scenario 3 achieved the highest coverage which addresses 8 out of 17 spill events. Conversely, Scenario 1 and Scenario 4 achieved moderate coverage 5 out of 17, whereas Scenario 5, despite receiving significant transfers, also reached 5 spills. The uncovered spills primarily occurred in remote or low sensitivity regions with marginal objective weight or were cost prohibitive under existing resource constraints.\\

Resource deployment patterns reveal clear distinctions in the operational role of each station. Station 3 served as the central deployment hub for most scenarios which leveraged its high-capacity skimmer and boomer inventories to serve both direct and transferred demand. Station 4 provided secondary coverage which often relies on interstation transfers from station 3.
All covered spill events received 100\% of the required demand for each resource type like skimmers, booms, and dispersants across all scenarios. In total, more than 2,351 units of booms, 5.6 skimmers and 4.0 units of dispersants were allocated on average per scenario. Transfers were critical in balancing supply across stations with substantial quantities, specifically 80,000 units of booms and 15 skimmers were moved between stations to address time dependent or high priority spills.\\

Resource utilization analysis shows that even under optimal coverage, overall system utilization remained low to moderate, which is consistent with robust design principles in emergency logistics. Skimmer utilization peaked at 50\% at station 4 in Scenario 5, while boomer and dispersant utilization remained below 7\% at all stations. This suggests the model maintains operational slack to ensure responsiveness under demand spikes, which is crucial in uncertain and high stakes Arctic conditions.

\subsection{Value of Stochastic Solution}
\vspace{0.5em}
Three benchmark models were evaluated to quantify the benefit of incorporating uncertainty. Those are as follows: RP (Recourse Problem), the full stochastic model solved across all five scenarios, EVP (Expected Value Problem), a deterministic model using mean values for all uncertain parameters, and EEV (Expected value of EVP solution), the performance of the deterministic solution when tested under the five stochastic scenarios. Table \ref{tbl:voi_metrics} summarizes the comparative performance of these approaches:

\begin{table}[width=.75\linewidth,cols=3,pos=h]
\caption{Value of Information Metrics}
\label{tbl:voi_metrics}
\begin{tabular*}{\tblwidth}{@{} L C C @{}}
\toprule
\textbf{Metric} & \textbf{Value} & \textbf{Relative Performance (\%)} \\
\midrule
RP (Stochastic)        & 1.9855 & 100.00 \\
EWS (Perfect Information) & 1.9855 & 100.00 \\
EVP (Deterministic)    & 1.7855 & 89.92  \\
EEV (Expected Value)   & 1.4659 & 73.83  \\
VSS                    & 0.5196 & 35.45  \\
EVPI                   & 0.0000 & 0.00   \\
\bottomrule
\end{tabular*}
\end{table}

The VSS defined as the difference between the RP and EEV is 0.5196 which equivalent to a 35.45\% gain in objective performance over deterministic planning. This substantial value demonstrates the benefit of explicitly incorporating scenario uncertainty in the model. In contrast, the Expected Value of Perfect Information (EVPI) is approximately zero which indicates that the stochastic solution performs equivalently to an ideal solution with full scenario knowledge. Therefore, no meaningful gains would be realized from perfect foresight which, thereafter, reinforces the adequacy and robustness of the stochastic model.

\subsection{Scenario Level Insights}
\vspace{0.5em}
Table \ref{tbl:scenario_objective_comparison} presents a comparison between the performance of the deterministic solution (EVP) and the ideal perfect-information model (EWS) under each scenario. The objective value for each approach is shown alongside the corresponding scenario probability.\\

\begin{table}[width=.75\linewidth,cols=4,pos=h]
\caption{Scenario-wise Objective Comparison}
\label{tbl:scenario_objective_comparison}
\begin{tabular*}{\tblwidth}{@{} C C C C @{}}
\toprule
\textbf{Scenario} & \textbf{EWS Value} & \textbf{EVP Evaluation} & \textbf{Probability} \\
\midrule
1 & 1.8611 & 1.2763 & 0.3751 \\
2 & 1.4641 & 1.2186 & 0.2178 \\
3 & 2.7776 & 0.6636 & 0.0696 \\
4 & 2.1935 & 2.1935 & 0.1980 \\
5 & 2.4440 & 1.7294 & 0.1395 \\
\bottomrule
\end{tabular*}
\end{table}

Scenario 3 illustrates the most significant divergence, where the deterministic solution captures only 23.9\% of the objective value achievable under perfect information. This underperformance stems from the inability of the deterministic model to anticipate high-magnitude, low-probability spill events, leading to under-allocation of critical resources. Conversely, Scenario 4 shows complete alignment between both models, indicating that in some cases, deterministic planning may suffice when variability is low, or response priorities are clear.\\

\begin{figure}[pos=H]
    \centering
    \includegraphics[width=0.70\textwidth,
    height=0.25\textheight,
    keepaspectratio]{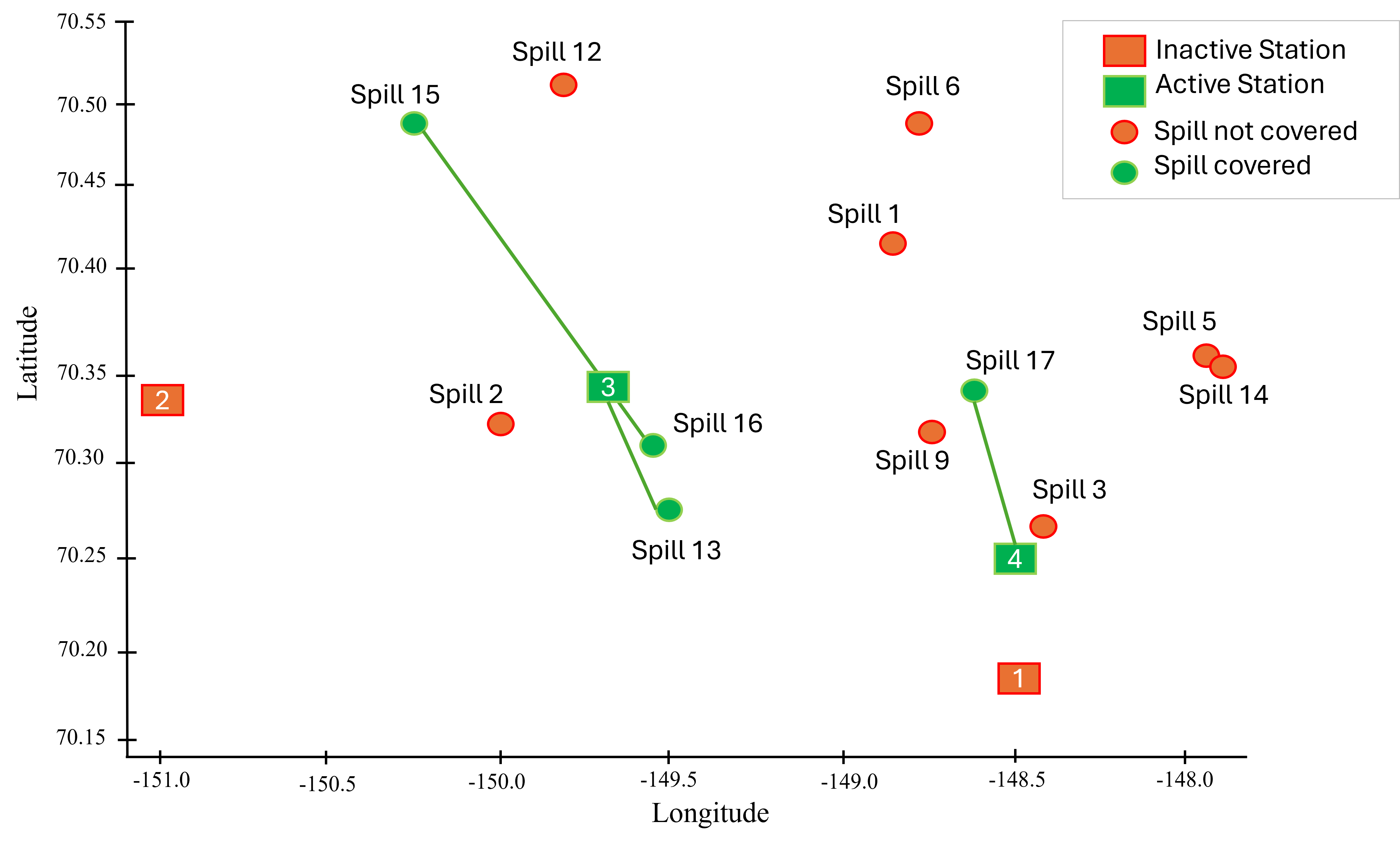}%
    \caption{Multi-resource Arctic oil spill response results for
    scenario~1}
    \label{fig:scenario1}
\end{figure}

\begin{figure}[pos=H]
    \centering
    \includegraphics[width=0.70\textwidth,
    height=0.25\textheight,
    keepaspectratio]{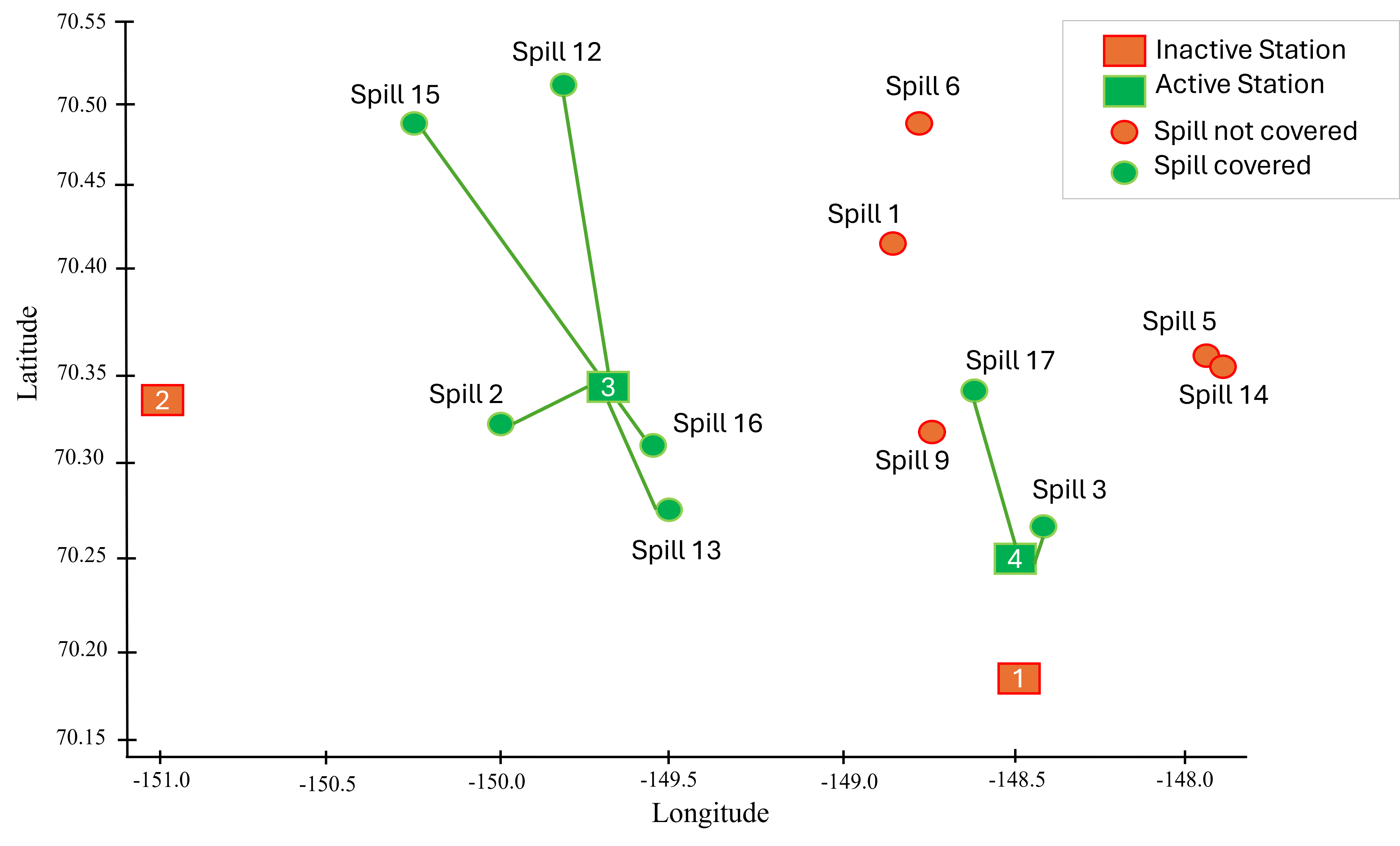}
    \caption{Multi-resource Arctic oil spill response results for
    scenario~2}
    \label{fig:scenario2}
\end{figure}

\begin{figure}[pos=H]
    \centering
    \includegraphics[width=0.70\textwidth,
    height=0.25\textheight,
    keepaspectratio]{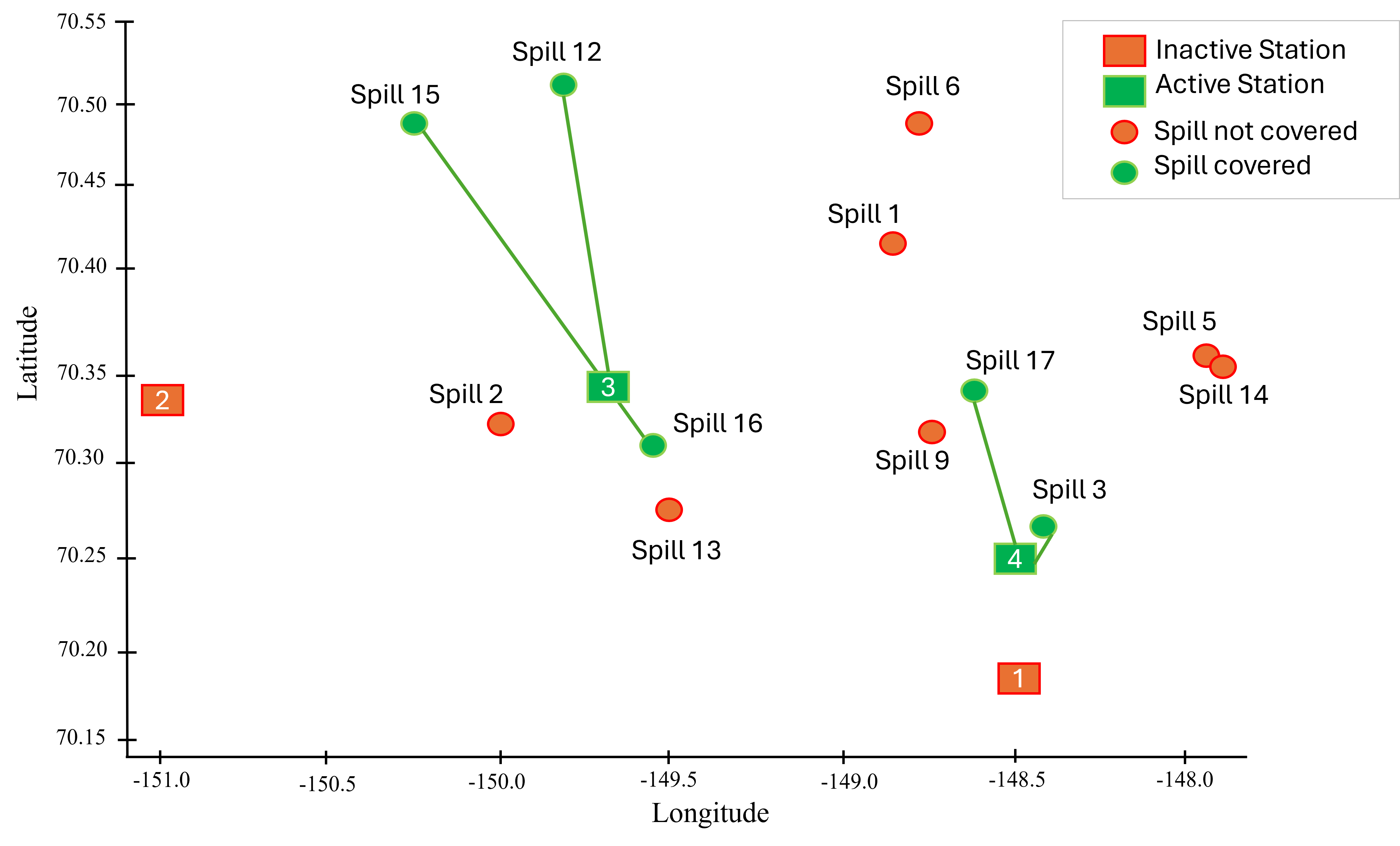}
    \caption{Multi-resource Arctic oil spill response results for
    scenario~3}
    \label{fig:scenario3}
\end{figure}

\begin{figure}[pos=H]
    \centering
    \includegraphics[width=0.70\textwidth,
    height=0.25\textheight,
    keepaspectratio]{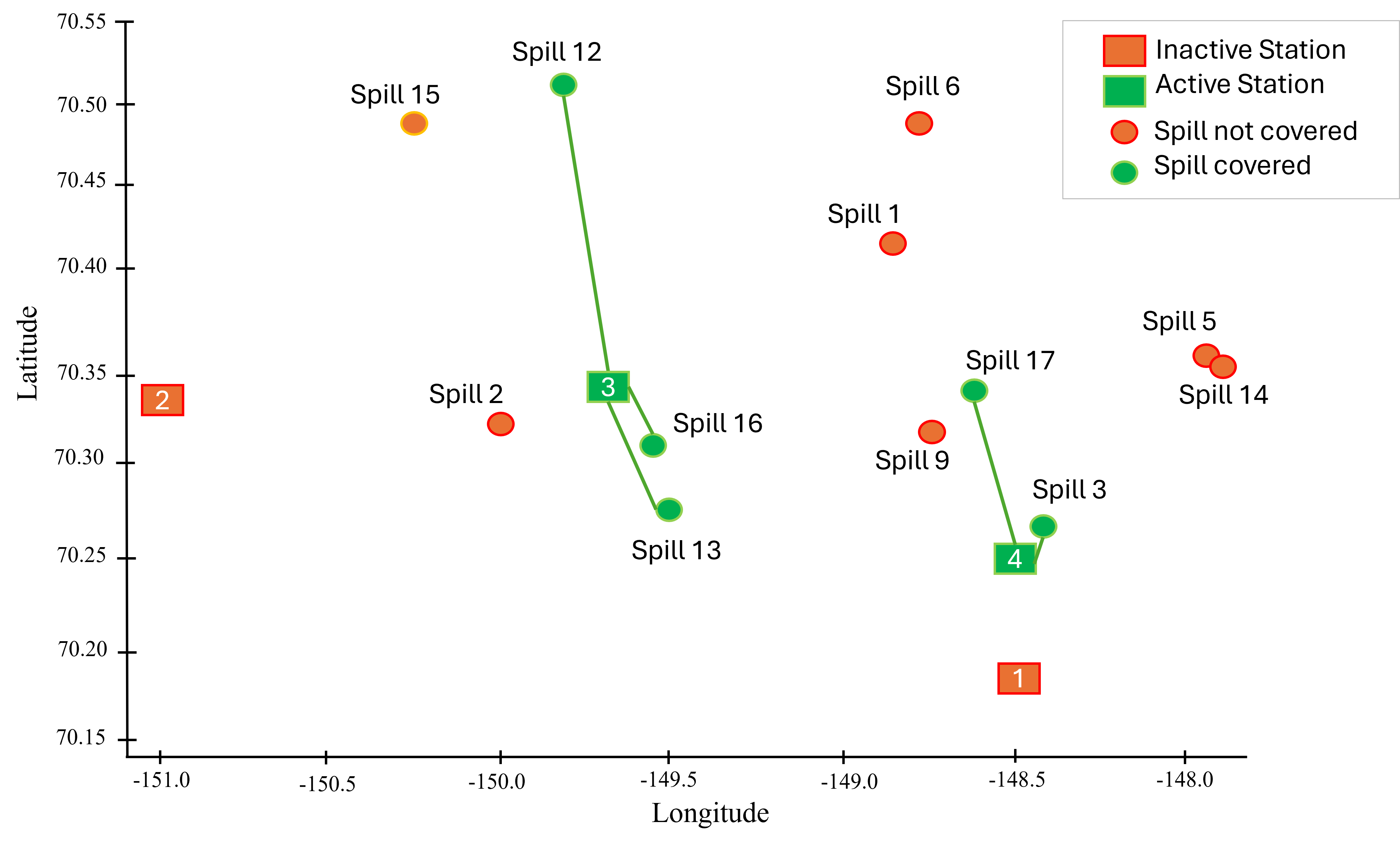}
    \caption{Multi-resource Arctic oil spill response results for
    scenario~4}
    \label{fig:scenario4}
\end{figure}

\begin{figure}[pos=H]
    \centering
    \includegraphics[width=0.70\textwidth,
    height=0.3\textheight,
    keepaspectratio]{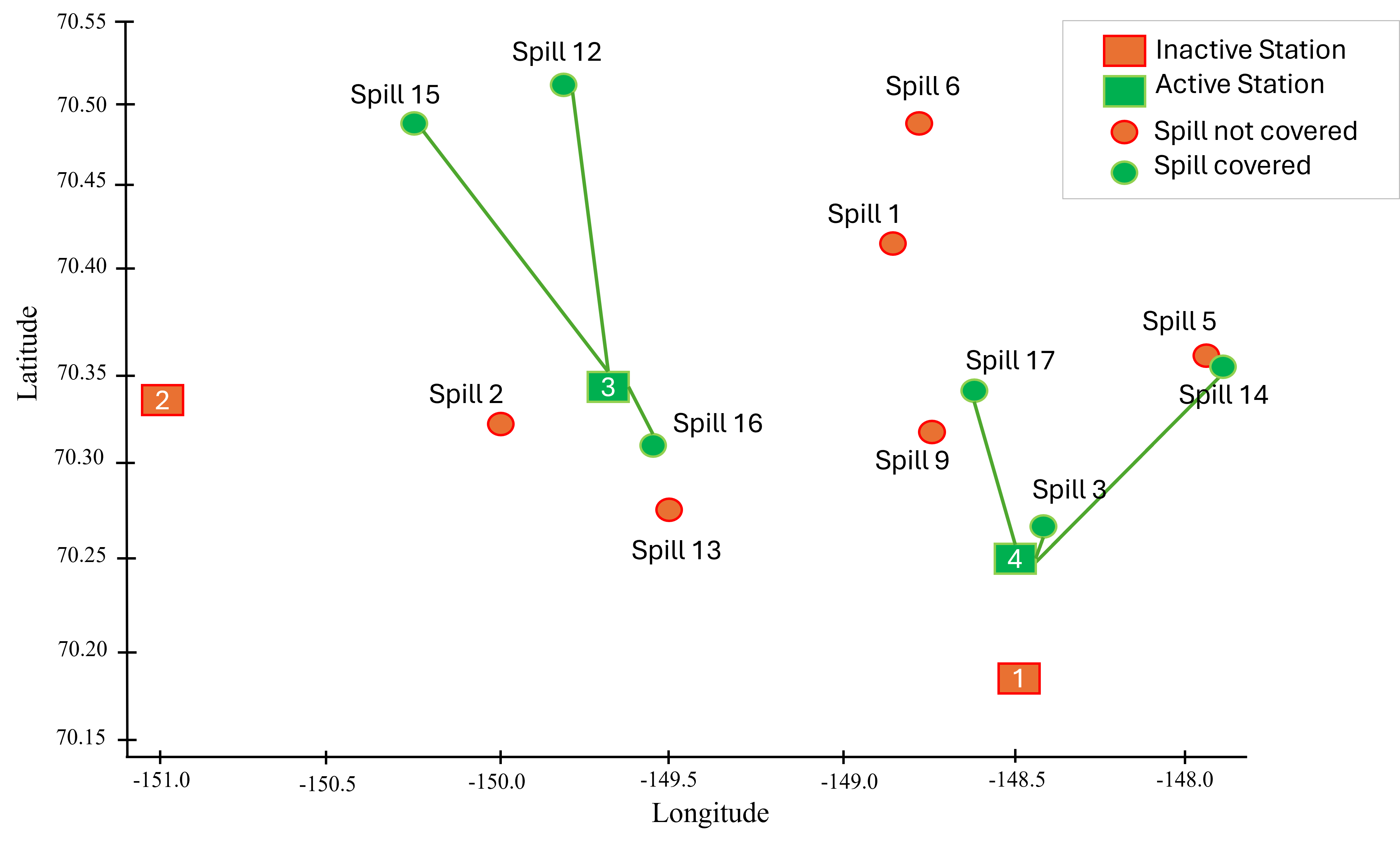}
    \caption{Multi-resource Arctic oil spill response results for
    scenario~5}
    \label{fig:scenario5}
\end{figure}

Figures~\ref{fig:scenario1}–\ref{fig:scenario5} demonstrate the spatial
assignments, and scenario-wise performance summaries produced by the multi-resource stochastic model. These visuals confirm that stochastic planning enables more targeted coverage of high-impact events, reduces idle or wasted capacity, and enhances overall operational resilience in the Arctic response network.\\

From a strategical standpoint, the findings emphasize the importance of deploying a flexible and risk informed response system in the Arctic. The results suggest that a two-station network, when optimized under uncertainty, can deliver high-impact coverage while maintaining sufficient spare capacity. Additionally, the negligible EVPI implies that while improved forecasting has limited marginal value, incorporating uncertainty within the decision model itself provides significant benefit. Therefore, investment should focus on robust infrastructure, logistical coordination, and real time deployment capabilities rather than predictive analytics alone.

\section{Sensitivity Analysis}
\vspace{0.5em}
A comprehensive sensitivity analysis was performed to evaluate the robustness and responsiveness of the proposed multi resource stochastic Arctic oil spill response model. This analysis involved 324 model runs across varying global objective weights of coverage importance ($k_1$) and cost importance ($k_2$) and internal weights assigned to the three primary spill response criteria including spill size ($\omega_1$), environmental sensitivity ($\omega_2$), and response time ($\omega_3$). The objective was to quantify how tradeoffs among competing factors influence station selection, resource allocation, spill coverage, and cost efficiency. All results were normalized to ensure comparability, and outputs were analyzed both statistically and visually to extract actionable insights.

\subsection{Effect of Coverage Importance Weight ($k_1$) on Coverage Objective Value}
\begin{figure}[pos=H]
	\centering
	\includegraphics[width=0.7\textwidth,
	height=0.25\textheight,
	keepaspectratio]{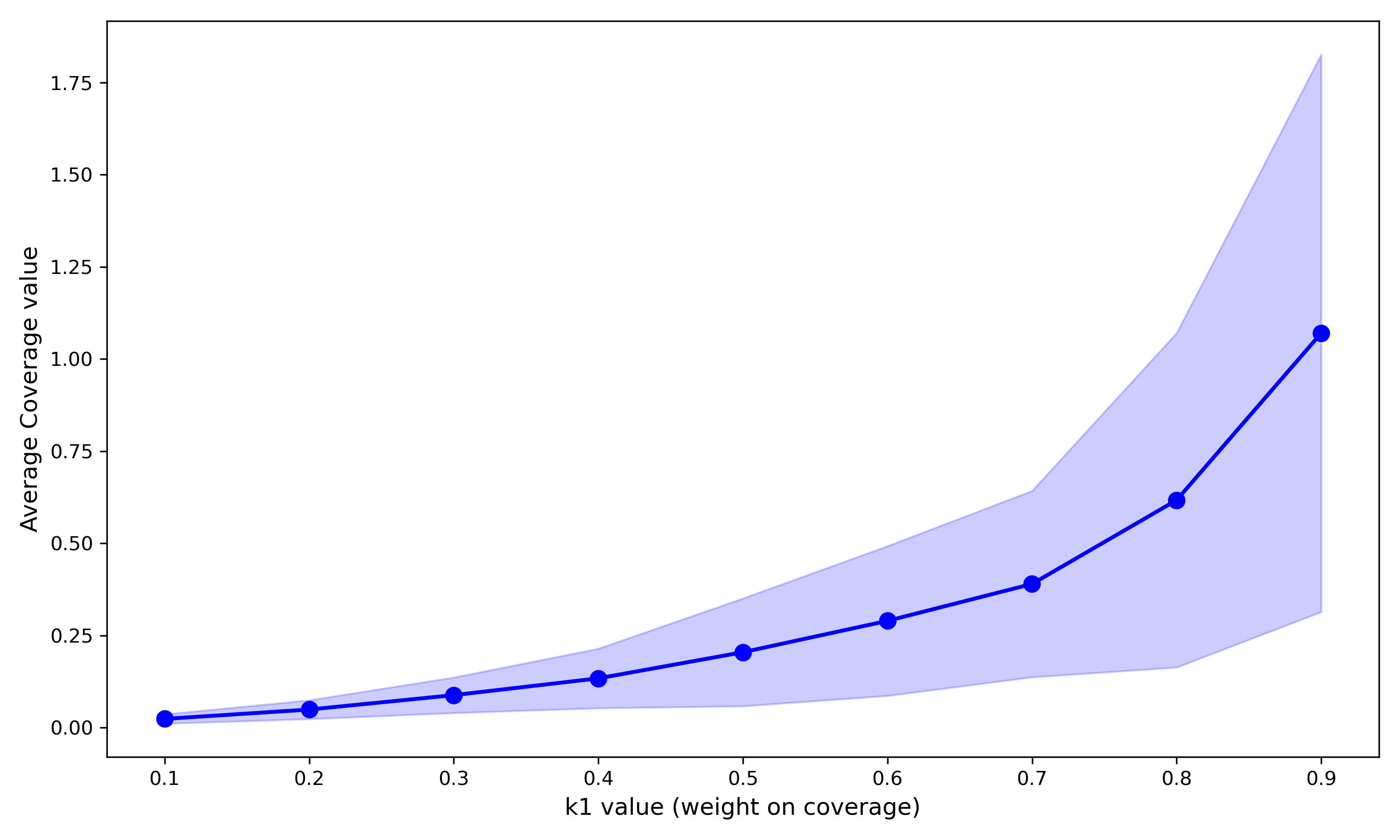}
	\caption{Effect of coverage importance weight ($k_1$) on coverage objective value.}
	\label{fig:k1coverage}
\end{figure}

\vspace{0.5em}
The parameter $k_1$ controls the global emphasis on maximizing spill coverage in the model's objective function. As illustrated in figure \ref{fig:k1coverage}, there is a significant positive correlation between $k_1$ and the average weighted coverage value. When $k_1$ is set to a low value, say between 0.1–0.3, the model demonstrates limited spill coverage, as the objective function is dominated by cost saving purposes. Under these settings, fewer stations are selected, and minimal deployment occurs due to the high penalty associated with station, deployment, and transfer costs.\\

However, as $k_1$ increases beyond 0.5, the model begins prioritizing spill coverage more aggressively. At $k_1$ = 0.7, the average coverage value increases significantly, and by $k_1$ = 0.9, the average weighted coverage reaches 1.08 with high variance. This variability suggests that as the weight on coverage increases, the solution becomes more sensitive to scenario specific parameters such as spill location, resource availability, and environmental sensitivity. The standard deviation bands shown in the shaded region confirm that while the upper bounds of performance improve significantly with $k_1$, the model also explores a broader solution space which results in performance dispersion.

\subsection{Effect of Coverage Importance Weight ($k_1$) on the Overall Objective Value}
\begin{figure}
	\centering
	\includegraphics[width=0.7\textwidth,
	height=0.25\textheight,
	keepaspectratio]{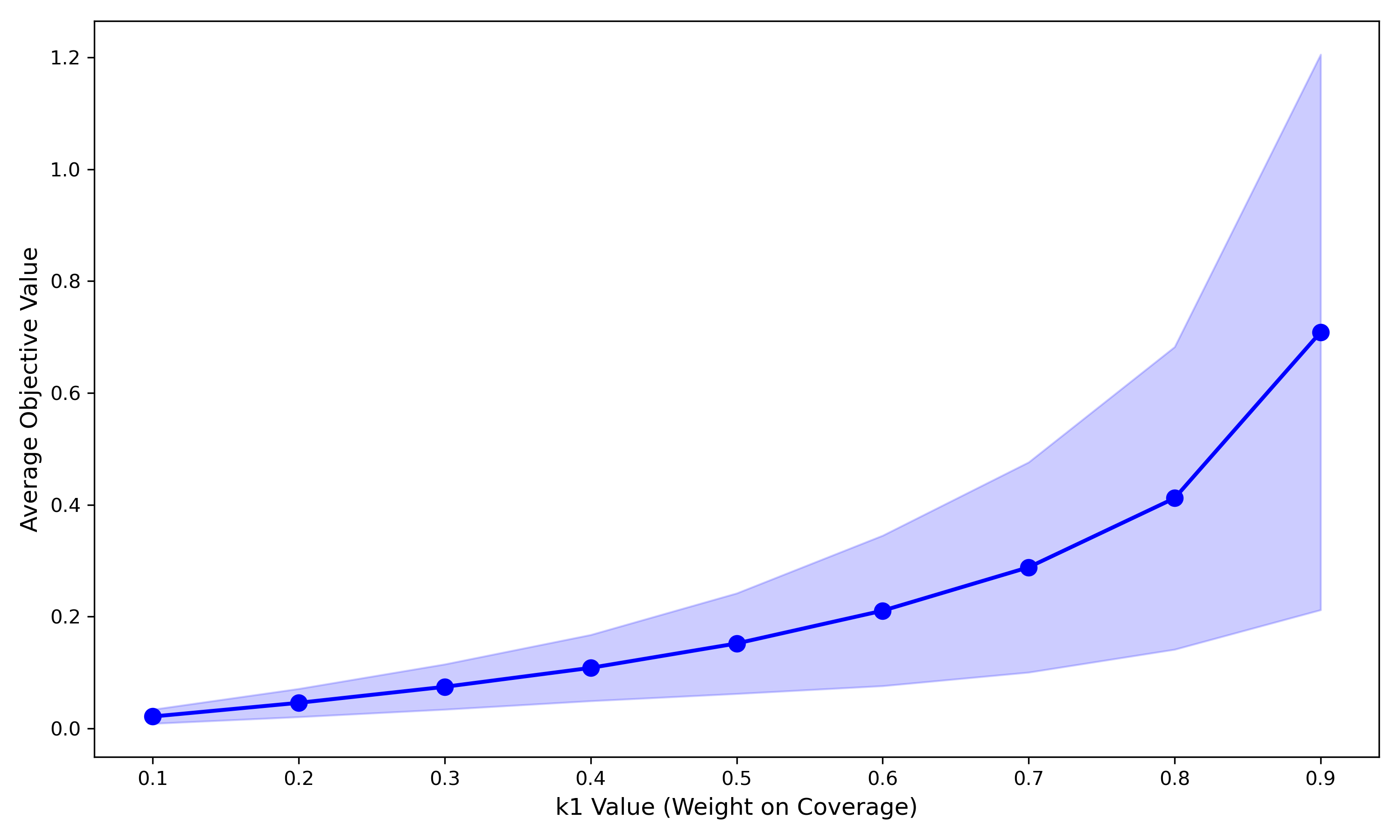}
	\caption{Effect of coverage importance weight ($k_1$) on overall objective.}
	\label{fig:k1coverageoverall}
\end{figure}

\vspace{0.5em}
Figure \ref{fig:k1coverageoverall} presents the effect of increasing $k_1$ on the overall objective value, which includes both weighted coverage and weighted cost terms. The results reinforce the previous trend. As $k_1$ increases, so does the average objective value. At $k_1$ = 0.1, the mean objective value remains below 0.05, which indicates negligible coverage and cost avoidance. From $k_1$ = 0.5 onwards, the objective value begins to rise more sharply, with the steepest increase observed between $k_1$ = 0.7 and 0.9. The highest mean objective value is 0.75 at $k_1$ = 0.9.\\

This indicates that assigning greater importance to coverage significantly enhances the system’s ability to meet spill response goals. Moreover, the increased standard deviation at high $k_1$ levels shows that such configurations enable a broader and more flexible range of decisions, particularly in selecting stations that may not be cost-efficient but offer strategic value for sensitive or high-volume spill response.

\subsection{Effect of Cost Importance Weight ($k_2$)}
\vspace{0.5em}
Since $k_2 = 1 - k_1$, a complementary trend emerges with respect to cost prioritization. When $k_2$ is high while $k_1$ is low, the model tends to suppress all high-cost operations, including opening stations and deploying resources which results in low or zero spill coverage. In multiple runs under $k_2 \geq 0.7$, the model achieved near zero objective values due to full avoidance of deployment costs, essentially ignoring the spills altogether.\\

\begin{figure}[pos=H]
	\centering
	\includegraphics[width=0.75\textwidth,
	height=0.3\textheight,
	keepaspectratio]{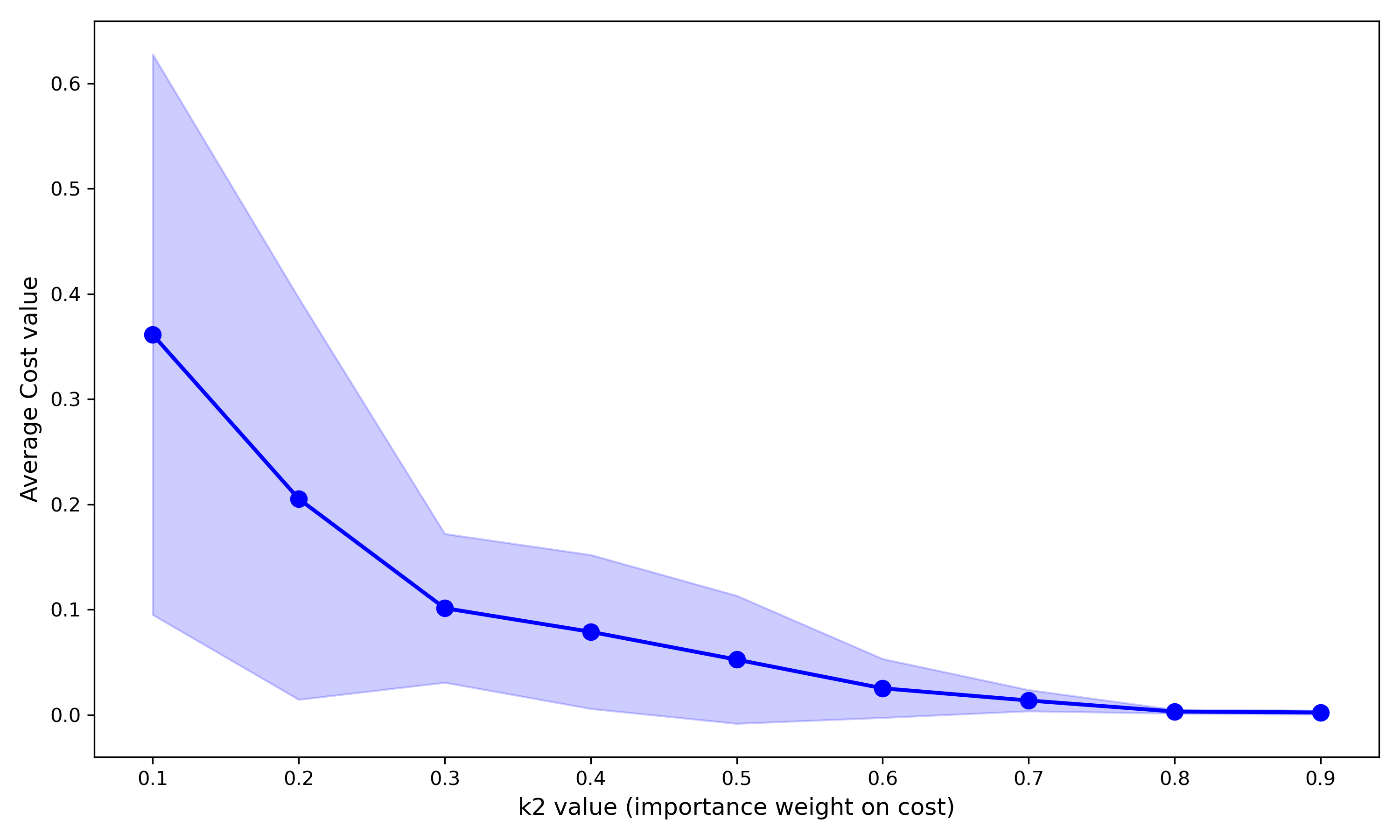}
	\caption{Effect of cost importance weight ($k_2$) on cost objective value.}
	\label{fig:k2cost}
\end{figure}

In contrast, low $k_2$ values ($\leq  0.2$) permitted more aggressive spill coverage actions, though at higher normalized cost values. This confirms that cost minimization and ecological performance are inherently conflicting objectives in this model. Without regulatory, ecological, or political mandates emphasizing coverage, solutions tend to gravitate toward low-cost but ineffective configurations. Figure \ref{fig:k2cost} shows the effect of cost importance on weight ($k_2$).

\subsection{Joint Sensitivity Analsysis: Coverage Internal Weight ($w_1$, $w_2$, $w_3$)}
\begin{figure}[pos=H]
    \centering
    \includegraphics[width=0.80\textwidth,
    height=0.30\textheight,
    keepaspectratio]{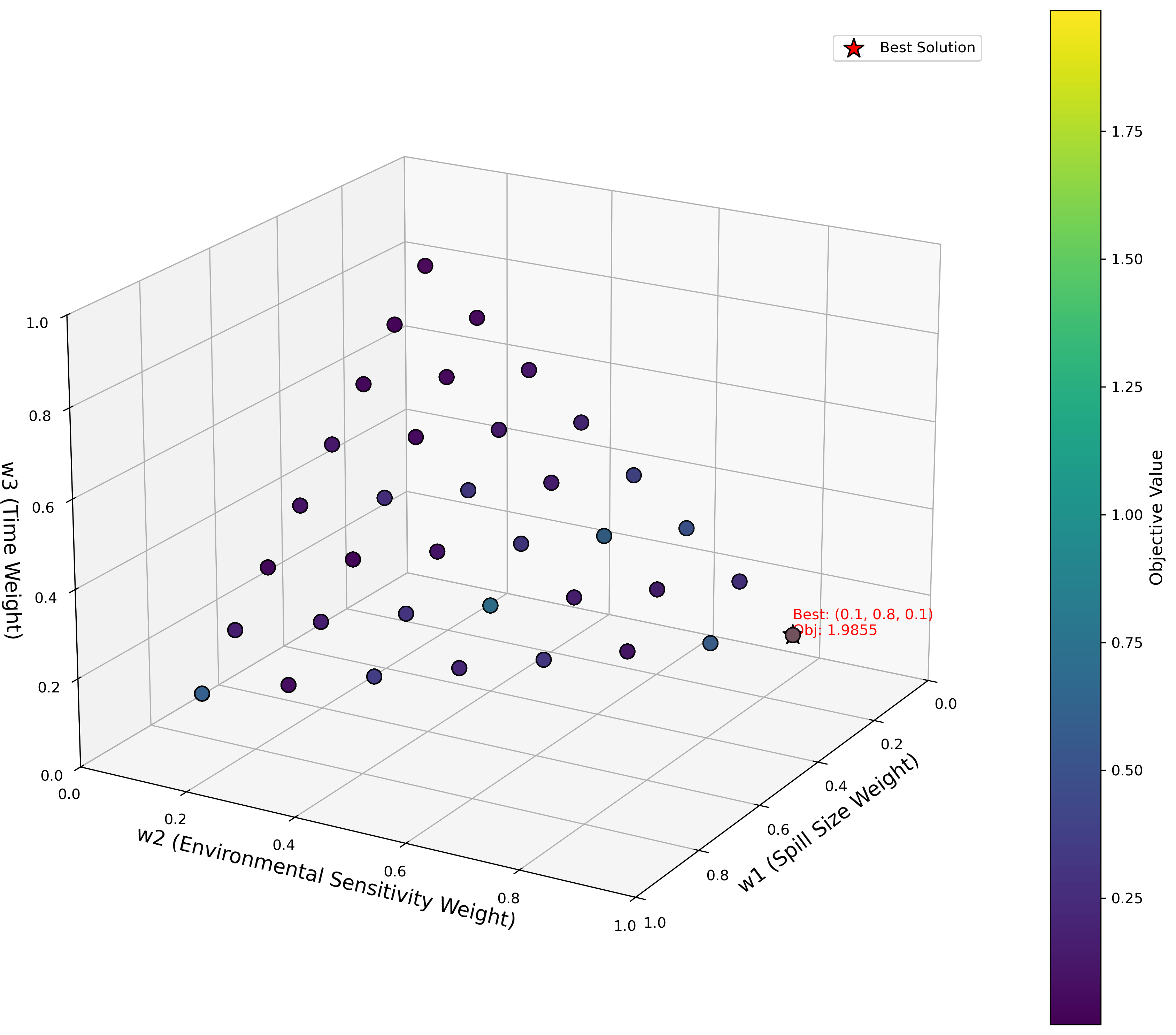}
    \caption{Joint Probability Analysis: 3D scatter visualization}
    \label{fig:jpweights}
\end{figure}

\vspace{0.5em}
A joint sensitivity analysis was performed across 81 valid combinations of the internal weights ($\omega_1$, $\omega_2$, $\omega_3$), given that $\omega_1 + \omega_2 + \omega_3 = 1$,  to better understand how individual components of spill coverage influence the overall optimization model performance.  These weights represent distinct environmental priorities: $\omega_1$ emphasizes spill size, giving preference to larger incidents; $\omega_2$ targets environmental sensitivity, prioritizing responses in ecologically fragile or high-risk zones; and $\omega_3$ corresponds to response time, prioritizing rapid deployment. The results, illustrated through a 3D scatter visualization, reveal that the optimal weight configuration was ($\omega_1$ = 0.1, $\omega_2$ = 0.8, $\omega_3$ = 0.1), which yielded the highest objective value of 1.9855. This configuration strongly indicates that environmental sensitivity is the most influential factor in designing effective spill response strategies in Arctic environments. Across all tested combinations, those with high $\omega_2$ and moderate $\omega_1$ values consistently outperformed others, while elevated values of $\omega_3$, especially beyond 0.3 tended to reduce solution quality. These results suggest that while rapid response is important, overemphasis on speed may compromise attention to environmental vulnerabilities and spill magnitude. Figure \ref{fig:jpweights} represents the joint sensitivity analysis through 3D scatter visualization.

\subsection{Pareto Frontier Analysis: Coverage vs. Cost}
\begin{figure}[pos=h]
	\centering
	\includegraphics[width=0.85\textwidth,
	height=0.4\textheight,
	keepaspectratio]{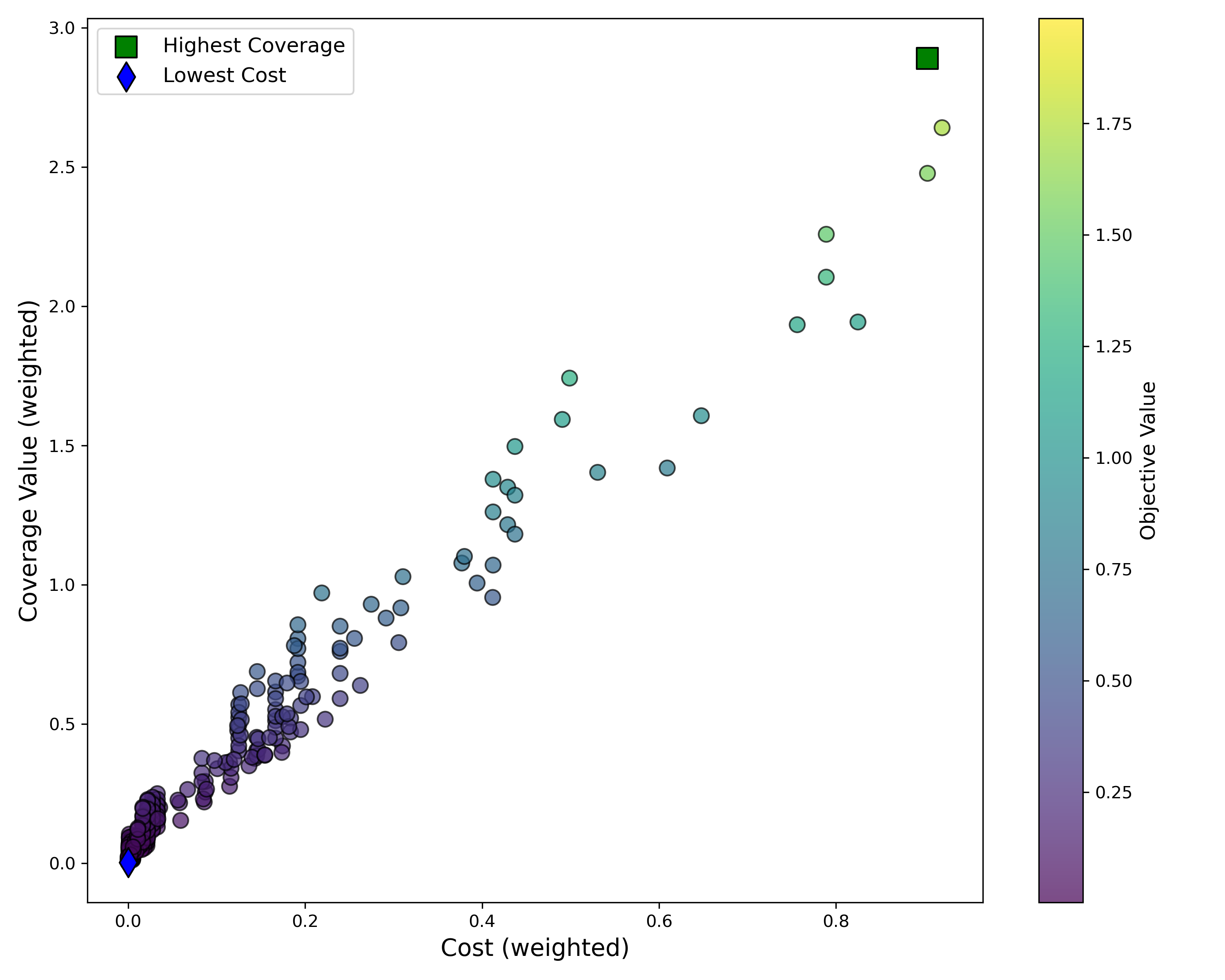}
	\caption{Pareto Frontier – Coverage vs Cost}
	\label{fig:pareto}
\end{figure}

\vspace{0.5em}
Figure \ref{fig:pareto} shows the Pareto frontier of all 324 solutions in the coverage cost space. Each point corresponds to a unique configuration of ($k_1$, $k_2$, $\omega_1$, $\omega_2$, $\omega_3$) which is colored by the resulting objective value. The frontier reveals a clear trade-off surface. Solutions on the lower left region minimize cost but offer poor coverage, whereas those in the upper right maximize coverage at significant expense.\\

The highest coverage solution with a weighted coverage value of 2.88 corresponds to the optimal parameter set $k_1$ = 0.9, $\omega_1$ = 0.1, $\omega_1$ = 0.8, $\omega_3$ = 0.1. On the other hand, the lowest cost solution appears near zero on the x-axis and exhibits negligible coverage. The shape of the pareto frontier demonstrates that for moderate increases in cost, especially between normalized values of 0.1 and 0.4, the model can achieve substantial coverage improvements which can make these mid-range solutions particularly attractive for budget constrained but environmentally conscious planners.\\

The sensitivity analysis revealed several critical insights regarding the interplay between objective priorities and model performance. First, adjusting the global objective weights showed that increasing emphasis on spill coverage ($k_1$ $\geq$ 0.8) significantly enhanced both the overall objective value and the system’s ability to protect against spill events. Conversely, placing a high weight on cost minimization ($k_2$ $\geq$ 0.7) led to environmentally suboptimal outcomes. Second, among the internal coverage weights, environmental sensitivity ($\omega_2$) proved to be the most influential factor, which outweigh the effects of both spill size ($\omega_1$) and response time ($\omega_3$) in determining optimal response strategies. Third, the Pareto frontier analysis revealed a nonlinear tradeoff between cost and coverage. While moderate increases in cost led to substantial gains in spill coverage, the rate of return diminished beyond certain thresholds. Finally, the optimal configuration was found at $k_1$ = 0.9, $\omega_1$ = 0.1, $\omega_2$ = 0.8, $\omega_3$ = 0.1, yielding the highest objective value of 1.9855 and a corresponding coverage value of 2.88. This indicates that high coverage, environmentally sensitive strategies are not only effective but also achievable with moderate cost increases. Overall, these findings strongly advocate for a coverage centric, environmentally weighted planning framework in Arctic oil spill response operations, particularly in decision-making environments that prioritize long-term ecological sustainability alongside economic considerations.

\subsection{Discussion}
\vspace{0.5em}
Based on the findings of this study, integrating stochastic optimization into Arctic oil spill response planning is strongly recommended. The model demonstrates a substantial 35.45\% performance increase over traditional deterministic approaches, which proves inadequate for low-probability, high-impact scenarios. Key operational insights confirm the model's practicality. The consistent selection of Stations 3 and 4 provides robust, complementary coverage. Furthermore, the analysis validates that low resource utilization rates are not inefficient but are essential reserve capacity for resilience. Strategically, the sensitivity analysis reveals that environmental sensitivity is the most critical factor, and over-prioritizing rapid deployment can be counterproductive. The Pareto frontier shows that mid-range solutions offer the best balance between ecological coverage and cost. Finally, since the stochastic solution nearly matches the performance of a perfect-information model, planners should prioritize risk-informed infrastructure and flexible logistics over investing in better scenario forecasting.\\

\subsection*{Acknowledgments}
\vspace{0.4em}
 Special thanks to Alaska Clean Seas (ACS) and NOAA’s Environmental Sensitivity Index (ESI) program, whose datasets and reports were essential in constructing the oil spill response model.

\subsection*{Data availability}
\vspace{0.4em}
Data will be made available upon request.

\subsection*{Funding}
\vspace{0.4em}
No dedicated funding was received for this research

\subsection*{Declaration of competing interest}
\vspace{0.4em}
The authors declare that they have no known competing financial interests or personal relationships that could have appeared to influence the work reported in this paper.

\subsection*{Credit authorship contribution statement}
\vspace{0.4em}
\textbf{Md Ashiqur Rahman}: Conceptualization, Methodology, Investigation, Software, Formal analysis, Validation, Data curation, Visualization, Writing – original draft, Writing – review \& editing, \textbf{Mustofa Tanbir Kuhel}:  Methodology, Investigation, Data curation, Visualization, Writing – original draft, Writing – review \& editing, \textbf{Dr. Clara M Novoa Ramirez}: Conceptualization, Methodology, Supervision, Data Curation, Formal Analysis, Validation, Writing – review \& editing.

\newpage
%% Keep Difference %%%
\bibliographystyle{cas-model2-names}
\bibliography{cas-refs}

\end{document}